\documentclass[twoside]{article}

\usepackage{amsmath,amsthm,amssymb}

\usepackage{newtxtext,newtxmath}

\usepackage[text={12.5cm,19cm},centering,paperwidth=17cm,paperheight=24cm]{geometry}

\usepackage[fontsize=10.4pt]{scrextend}

\def\titlerunning#1{\gdef\titrun{#1}}

\makeatletter
\def\author#1{\gdef\autrun{\def\and{\unskip, }#1}\gdef\@author{#1}}
\def\address#1{{\def\and{\\\hspace*{15.6pt}}\renewcommand{\thefootnote}{}\footnote{#1}}\markboth{\autrun}{\titrun}}
\makeatother

\def\email#1{email: \href{mailto:#1}{#1} }
\def\subjclass#1{\par\bigskip\noindent\textbf{Mathematics Subject Classification 2020.} #1}
\def\keywords#1{\par\smallskip\noindent\textbf{Keywords.} #1}

\newenvironment{acknowledgments}{\bigskip\small\noindent\textit{Acknowledgments.}}{\par}

\newtheorem{thm}{Theorem}[section]

\newtheorem{conjecture}{Conjecture}

\theoremstyle{definition}

\newtheorem{exa}[thm]{Example}

\newtheorem*{rem}{Remark}

\numberwithin{equation}{section}

\usepackage[hyperfootnotes=false,colorlinks=true,allcolors=blue]{hyperref}

\begin{document}

\titlerunning{Geometry of Spherical Harmonics}

\title{\textbf{Some Recent Developments on the Geometry of Random Spherical
Eigenfunctions}}

\author{Domenico Marinucci}

\date{}

\maketitle

\address{Department of Mathematics, University of Rome Tor Vergata, \email{marinucc@mat.uniroma2.it}}


\begin{abstract}

A lot of efforts have been devoted in the last decade to the investigation
of the high-frequency behaviour of geometric functionals for the excursion
sets of random spherical harmonics, i.e., Gaussian eigenfunctions for the
spherical Laplacian $\Delta_{\mathbf{S}^2}$. In this survey we shall review some of these results, with particular
reference to the asymptotic behaviour of variances, phase transitions in the
nodal case (the \emph{Berry's Cancellation Phenomenon}), the distribution of the
fluctuations around the expected values, and the asymptotic correlation
among different functionals. We shall also discuss some connections with the
Gaussian Kinematic Formula, with Wiener-Chaos expansions and with recent
developments in the derivation of Quantitative Central Limit Theorems (the
so-called Stein-Malliavin approach).

\end{abstract}

\section{Introduction}

Spherical eigenfunctions are defined as the solutions of the Helmholtz
equation
\begin{equation*}
\Delta _{\mathbb{S}^{2}}f_{\ell }+\lambda _{\ell }f_{\ell }=0\text{ },\text{
}f_{\ell }:\mathbb{S}^{2}\rightarrow \mathbb{R}\text{ },\text{ }\ell
=1,2,...,
\end{equation*}%
where $\Delta _{\mathbb{S}^{2}}$ is the spherical Laplacian and $\left\{
-\lambda _{\ell }=-\ell (\ell +1)\right\} _{\ell =1,2,...}$ is the set of
its eigenvalues. A random structure can be constructed easily by assuming
that the eigenfunctions $\left\{ f_{\ell }(\cdot)\right\} $ follow a Gaussian
isotropic random process on $\mathbb{S}^{2}.$ More precisely, for each $x\in
\mathbb{S}^{2}$, we take $f_{\ell }(x)$ to be a Gaussian
random variable defined on a suitable probability space $\left\{ \Omega ,\Im
,\mathbb{P}\right\} $; without loss of generality, we assume $\left\{
f_{\ell }(\cdot)\right\} $ to have mean zero, unit variance, and covariance
function given by
\begin{equation*}
\mathbb{E}\left[ f_{\ell }(x)f_{\ell }(y)\right] =P_{\ell }(\left\langle
x,y\right\rangle) \text{ , }x,y\in \mathbb{S}^{2},\text{ }P_{\ell }(t):=%
\frac{1}{2^{\ell }\ell !}\frac{d^{\ell }}{dt^{\ell }}(t^{2}-1)\text{ , } t \in \left[ -1,1 \right ] \text{ ,}
\end{equation*}%
where $\left\{ P_{\ell }(.)\right\} $ denotes the family of Legendre
polynomials: this is the only covariance structure to ensure that the random
eigenfunctions are isotropic, that is, invariant in law with respect to the
action of the group of rotations $SO(3)$. Random spherical eigenfunctions,
also known as random spherical harmonics, arise in a huge number of
applications, especially in connection with Mathematical Physics: in
particular, their role in Quantum Chaos has drawn strong interest in the
last two decades, starting from the seminal papers by \cite{Berry 1977},
\cite{Berry 2002}, \cite{Wig}, \cite{nazarov}; also, they represent the
Fourier components of isotropic spherical random fields, whose analysis
has an extremely important role in Cosmology (see e.g., \cite{MaPeCUP}). Of course, random spherical harmonics are just a
special case of a much richer literature on random eigenfunctions on general
manifolds; special interest has been
drawn for instance by \emph{Arithmetic Random Waves, }i.e., random
eigenfunctions on the torus $\mathbb{T}^{d},$ which were introduced by \cite%
{RW08} and then studied among others by \cite{KKW}, \cite{KW}, \cite%
{MPRW2015}, \cite{dalmao}, \cite{Rudnick}, \cite{RudnickYesha}, \cite{Cam19}%
, \cite{Maff}, \cite{Buckley}, see also \cite{CH}, \cite{Sarnak} and the
references therein. Although some of the results that we shall discuss have
related counterparts on the torus, on the higher-dimensional spheres, on more general compact manifolds and in the
Euclidean case, we will stick mainly to $%
\mathbb{S}^{2}$ for brevity and simplicity.

A lot of efforts have been spent in the last decade to characterize the
geometry of the excursion sets of random spherical harmonics, which are defined as
\begin{equation}
A_{u}(f_{\ell };\mathbb{S}^{2}):=\left\{ x\in \mathbb{S}^{2}:f_{\ell
}(x)\geq u\right\} \text{ , }u\in \mathbb{R}\text{ .}  \label{Excursion}
\end{equation}%
A classical tool for the investigation of these sets is given by the
so-called Lipschitz-Killing Curvatures (or equivalently, by Minkowski
functionals, see \cite{adlertaylor}), which in dimension 2 correspond to the
Euler-Poincar\'{e} characteristic, (half of) the boundary length and the
excursion area. A general expression for their expected values (covering
much more general Gaussian fields than random eigenfunctions) is given by
the \emph{Gaussian Kinematic Formula} (see \cite{Taylor},\cite{adlertaylor}%
). Over the last decade, more refined characterizations for random spherical
harmonics have been obtained, including neat analytic expressions (in the
high energy limit $\lambda _{\ell }\rightarrow \infty $) for the
fluctuations around their expected values and the correlation among these
different functionals; much of the literature has been concerned with the
\emph{nodal} case, corresponding to $u=0,$ to which we shall devote special
attention. In this survey, we shall review some of these results and present
some open issues for future research.

\section{The Gaussian Kinematic Formula for Lipschitz-Killing Curvatures on
Excursions Sets}

\subsection{The Kac-Rice Formula and the Expectation Metatheorem}

The first modern attempt to investigate the geometry of random processes and
fields can probably be traced back to the groundbreaking work by Kac (1943)
and Rice (1945) (\cite{Kac}, \cite{Rice}) on the zeroes of stochastic
processes. Their pioneering argument can be introduced as follows: let $%
f(.,.):\Omega \times \mathbb{R}\rightarrow \mathbb{R}$ be a continuous
stochastic process satisfying regularity conditions; our aim is to derive
the expected cardinality of its zero set in some finite interval (say $[0,T]$), i.e. the mean of%
\begin{equation*}
N_{0}([0,T]):=Card\left\{ t\in \mathbb{[}0,T\mathbb{]}:f(t)=0\right\} \text{
.}
\end{equation*}%
Now assume that $\left\{ f(.)\right\} $ is $C^{1}$ with probability one,
such that $f(0),f(T)\neq 0$ and
\begin{equation*}
\left\{ t:f(t)=0,\text{ }f^{\prime }(t)=0\right\} =\varnothing \text{ ;}
\end{equation*}%
then the following result (\emph{Kac's counting Lemma}) can be established
easily (see \cite{azaiswschebor}, p.69):%
\begin{equation*}
N_{0}([0,T])=\lim_{\varepsilon \rightarrow 0}\int_{0}^{T}\frac{1}{%
2\varepsilon }\mathbb{I}_{(-\varepsilon ,\varepsilon )}(f(t))|f^{\prime
}(t)|dt\text{ ,}
\end{equation*}%
where as usual $\mathbb{I}_{A}$ denotes the indicator function of the set $A$. With further efforts and assuming all exchanges of integrals and limits can
be justified, one obtains also%
\begin{equation}
\mathbb{E}\left[ N_{0}([0,T])\right] =\int_{0}^{T}\mathbb{E}\left[ \left.
|f^{\prime }(t)|\right\vert f(t)=0\right] p_{f(t)}(0)dt\text{ ,}
\label{KacRice}
\end{equation}%
where $\mathbb{E}\left[ \left. .\right\vert .\right] $ denotes as usual
conditional expected value and $p_{f}(.)$ the marginal density of $f(.),$
which is assumed to exist and admit enough regularity conditions (in the
overwhelming majority of the literature and in this whole survey, $f(.)$
will indeed be assumed to be Gaussian); (\ref{KacRice}) is the simplest
example of the \emph{Kac-Rice Formula.}

The basic idea behind the Kac-Rice approach has proved to be extremely
fruitful, leading to an enormous amount of applications and generalizations.
In particular, in the research monographs \cite{adlertaylor}, \cite%
{azaiswschebor}, (slightly different) versions of a general \emph{%
Expectation Metatheorem} (in the terminology of \cite{adlertaylor}) are
proved. More precisely, let us take $M$ to be a compact, $d-$dimensional
oriented $C^{1}$ manifold with a $C^{1}$ Riemannian metric $g.$ Assume $%
f:M\rightarrow \mathbb{R}^{d}$ and $h:M\rightarrow \mathbb{R}^{k}$ are
vector-valued random fields which satisfy suitable regularity conditions
(see \cite{adlertaylor}, \cite{azaiswschebor} for more details and \cite%
{Stecconi2021} for some very recent developments). Let $B\subset \mathbb{R}%
^{k}$ be a subset with boundary dimension smaller or equal than $k-1;$ then
define%
\begin{equation*}
N_{u}(f,h,M,B)=\left\{ t\in M:f(t)=u,h(t)\in B\right\} \text{ , }u\in
\mathbb{R}^{d}.
\end{equation*}%
The following extension of the Kac-Rice formula holds:

\begin{thm}
\label{Metatheorem} (\cite{adlertaylor}, \cite{azaiswschebor}) We have that%
\begin{equation*}
\mathbb{E}\left[ N_{u}(f,h,M,B)\right] =\int_{M}\mathbb{E}\left[ \left.
\left\vert \det \left\{ \nabla f(t)\right\} \right\vert \mathbb{I}%
_{B}(h(t))\right\vert f(t)=u\right] p_{f(t)}(u)\sigma _{g}(dt)\text{ ,}
\end{equation*}%
where as before $\mathbb{I}_{B}(.)$ denotes the indicator function, $\nabla f(.)$ the
(covariant) gradient of $f(.)$ and $\sigma _{g}(.)$ the volume form induced
by the metric $g$.
\end{thm}

\begin{rem}
By taking $k=1$, $f:=\nabla h$ the gradient of $h$ (and hence $\nabla
f=\nabla ^{2}h$ its Hessian) and $u=(0,...,0),$ Theorem \ref{Metatheorem}
yields the expected number of critical points with values in $B$ for the scalar random field $h$.
Simple modifications similarly yield the expected values for maxima, minima
and saddle points.
\end{rem}

The previous results have all been restricted to vector-valued random fields
whose image space has co-dimension zero. However, the results can be
similarly generalized to strictly positive co-dimensions. Indeed, under the
same setting as before assume instead that $f:M\rightarrow \mathbb{R}%
^{d^{\prime }}$ is such that $d^{\prime }<d;$ then $\nabla X$ is a $d\times
d^{\prime }$ rectangular matrix, and the following generalization of the
Expectation Metatheorem holds (see \cite{adlertaylor}, \cite{azaiswschebor})

\begin{thm}
\label{Metatheorem2} (\cite{adlertaylor}, \cite{azaiswschebor}) It holds that%
\begin{equation*}
\mathbb{E}\left[ \mathcal{H}_{u}(f,h,M,B)\right]
\end{equation*}%
\[
 =\int_{M}\mathbb{E}\left[
\left. \left\vert \det \left\{ (\nabla f(t))^{T}(\nabla f(t))\right\}
\right\vert ^{1/2}\mathbb{I}_{B}(h)\right\vert f(t)=u\right] p_{f(t)}(u)\sigma
_{g}(dt)\text{ ,}
\]
where $\mathcal{H}_{u}(f,h,M,B)$ denotes the $d-d^{\prime }$ dimensional
Hausdorff measure of the set $\left\{ t\in M:f(t)=u\text{ and }h(t)\in
B\right\} .$
\end{thm}

\begin{exa}
Let $M=\mathbb{S}^{2}$ the standard unit-dimensional sphere in $\mathbb{R}%
^{3},$ $f:\mathbb{S}^{2}\times \Omega \rightarrow \mathbb{R}$ a random
field, and let
\begin{equation*}
Len(f):=\mathcal{H}_{0}(f,\mathbb{S}^{2},0)=meas\left\{ t\in \mathbb{S}%
^{2}:f(t)=0\right\} \text{ ,}
\end{equation*}%
i.e., the length of the nodal lines of $f(.).$ Then
\begin{eqnarray*}
\mathbb{E}\left[ Len(f)\right] &=&\int_{\mathbb{S}^{2}}\mathbb{E}\left[
\left. \left\vert \det \left\{ (\nabla f(t))^{T}(\nabla f(t))\right\}
\right\vert ^{1/2}\right\vert f(t)=0\right] p_{f(t)}(0)\sigma (dt) \\
&=&\int_{\mathbb{S}^{2}}\mathbb{E}\left[ \left. \left\Vert \nabla
f(t)\right\Vert \right\vert f(t)=0\right] p_{f(t)}(0)\sigma (dt)\text{ ,}
\end{eqnarray*}%
where $\left\Vert .\right\Vert $ denotes Euclidean norm and $\sigma (.)$ the
standard Lebesgue measure on the unit sphere. In particular, assuming that
the law of $f(.)$ is isotropic (that is, invariant with respect to the
action of the group of rotations $SO(3))$ we get%
\begin{equation*}
\mathbb{E}\left[ Len(f)\right] =4\pi \times \mathbb{E}\left[ \left.
\left\Vert \nabla f(t)\right\Vert \right\vert f(t)=0\right] p_{f(t)}(0)\text{ .}
\end{equation*}
\end{exa}

\subsection{Intrinsic Volumes and Lipschitz-Killing Curvatures}

In the sequel, as mentioned earlier we will restrict our attention only to
Gaussian processes, which have driven the vast majority of research in this
area. We need now to introduce the Gaussian Kinematic Formula (see \cite%
{Taylor} and \cite{adlertaylor}); to this aim, let us first recall the
notion of \emph{Lipschitz-Killing Curvatures}. In the simplest setting of
convex subsets of the Euclidean space $\mathbb{R}^{d},$ Lipschitz-Killing
Curvatures (also known as intrinsic volumes) can be defined implicitly by
means of \emph{Steiner's Tube Formula}; to recall the latter, for any convex
set $d$-dimensional set $A\subset \mathbb{R}^{d}$ define the Tube of radius $%
\rho $ around $A$ as
\begin{equation*}
Tube(A,\rho ):=\left\{ x\in \mathbb{R}^{d}:d(x,A)\leq \rho \right\} \text{ ,
}d(x,A)=\inf_{y\in A}d(x,y)\text{ , }
\end{equation*}%
where $d(.,.)$ is the standard Euclidean distance. Then the following
expansion holds:%
\begin{equation*}
\mu _{d}\left\{ Tube(A,\rho )\right\} =\sum_{j=0}^{d}\omega _{d-j}\rho ^{d-j}%
\mathcal{L}_{j}(A)\text{ ,}
\end{equation*}%
where $\mathcal{L}_{j}(A)$ denotes the $j$-th Lipschitz-Killing Curvatures, $%
\mu _{d}(.)$ denotes the $d$-dimensional Lebesgue measure and $\omega _{j}:=%
\frac{\pi ^{j/2}}{\Gamma (\frac{j}{2}+1)}$ is the volume of the $j$%
-dimensional unit ball ($\omega _{0}=1,$ $\omega _{1}=2,$ $\omega _{2}=\pi ,$
$\omega _{3}=\frac{4}{3}\pi ).$

Lipschitz-Killing Curvatures can be shown to be additive and to scale with
dimensionality, in the sense that%
\begin{equation*}
\mathcal{L}_{j}(\lambda A)=\lambda ^{j}\mathcal{L}_{j}(A)\text{ for all }%
\lambda >0\text{ ,}
\end{equation*}%
and
\begin{equation*}
\mathcal{L}_{j}(A_{1}\cup A_{2})=\mathcal{L}_{j}(A_{1})+\mathcal{L}%
_{j}(A_{2})-\mathcal{L}_{j}(A_{1}\cap A_{2})\text{ .}
\end{equation*}%
For $j=d,$ it is immediately seen that $\mathcal{L}_{d}(A)$ is just the
Hausdorff measure of $A,$ whereas for $j=0$ we obtain $\mathcal{L}%
_{0}(A)=\varphi (A),$ the (integer-valued) Euler-Poincar\'e characteristic
of $A.$ A more general definition of $\mathcal{L}_{j}(.)$ can be given for
basic complexes (i.e., disjoint union of complex sets), for which the
following characterization (due to Hadwiger, see \cite{adlertaylor}) holds:%
\begin{equation}
\mathcal{L}_{j}(A)=\frac{\omega _{d}}{\omega _{d-j}\omega _{j}}\binom{d}{j}%
\int_{\mathcal{G}_{d}}\varphi (A\cap gE_{d-j})\mu (dg)  \label{Hadwiger}
\end{equation}%
where $\mathcal{G}_{d}=\mathbb{R}^{d}\times O(n)$ is the group of rigid motions, $%
E_{d-j}$ is any $d-j$ dimensional affine subspace and the volume form $\mu
(dg)$ is normalized so that
\begin{equation*}
\text{for all }x\in \mathbb{R}^{d},\text{ }A\subset \mathbb{R}^{d}\text{ , }%
\mu \left\{ g:gx\in A\right\} =\mathcal{H}(A)\text{ ,}
\end{equation*}%
where as before $\mathcal{H}(.)$ denotes the Hausdorff measure. For
instance, for $A=\mathbb{S}^{2}$ it is well-known and easy to check that (%
\ref{Hadwiger}) gives
\begin{equation*}
\mathcal{L}_{0}(\mathbb{S}^{2})=2\text{ , }\mathcal{L}_{1}(\mathbb{S}^{2})=0%
\text{ , }\mathcal{L}_{2}(\mathbb{S}^{2})=4\pi \text{ ,}
\end{equation*}%
which represent, respectively, the Euler-Poincar\'e characteristic, (half)
the boundary length and the area of the 2-dimensional unit sphere.

\subsection{The Gaussian Kinematic Formula}

From now on, we shall restrict our attention to Gaussian processes $%
f:M\rightarrow \mathbb{R}$, which we shall take to be zero-mean and
isotropic, meaning as usual that $\mathbb{E}\left[ f(t)\right] =0$ and $f(gt)%
\overset{d}{=}f(t)$ for all $t\in M\subset \mathbb{R}^{d}$ and $g\in SO(d);$
more explicitly, the law of the field $f(\cdot)$ will always be taken to be
invariant to rotations. In order to present the Gaussian Kinematic Formula,
let us first introduce a Riemannian structure governed by the covariance
function of the field $\left\{ f(.)\right\} $; more precisely, consider the
metric induced on the tangent plan $T_{t}M$ by the following inner product (%
\cite{adlertaylor}, p.305):
\begin{equation*}
g^{f}(X_{t},Y_{t}):=\mathbb{E}\left[ X_{t}f\cdot Y_{t}f\right] \text{ , }%
X_{t},Y_{t}\in T_{t}M\text{ . }
\end{equation*}%
This metric takes a particular simple form in case the field $f(.)$ is
isotropic; in these circumstances, $g^{f}(.,.)$ is simply the standard
Euclidean metric, rescaled by a factor that corresponds to the square root
of (minus) the derivative of the covariance density at the origin.

\begin{exa}
Consider the random spherical eigenfunction satisfying%
\begin{equation*}
\Delta f_{\ell }=-\lambda_{\ell}f_{\ell }\text{ , }f_{\ell }:\mathbb{S}%
^{2}\rightarrow \mathbb{R}\text{ , }\ell=0,1,2,...,
\end{equation*}%
with%
\begin{equation*}
\mathbb{E}\left[ f_{\ell }(x)\right] =0\text{ , }\mathbb{E}\left[ f_{\ell
}(x_{1})f_{\ell }(x_{2})\right] =P_{\ell }\left( \left\langle
x_{1},x_{2}\right\rangle \right) \text{ , }P_{\ell }^{\prime }(1)=-\frac{%
\ell (\ell +1)}{2}\text{ .}
\end{equation*}%
Then the induced inner product is simply%
\begin{equation*}
g^{f_{\ell }}(X,Y)=\sqrt{\frac{\ell (\ell +1)}{2}}\left\langle
X,Y\right\rangle _{\mathbb{R}^{3}}\text{ ;}
\end{equation*}%
this change of metric can of course be realized by transforming $\mathbb{S}%
^{2}$ into $\mathbb{S}_{\sqrt{\lambda _{\ell }/2}}^{2}:=\sqrt{\lambda _{\ell
}/2}\mathbb{S}^{2}.$
\end{exa}

Let us now write $\mathcal{L}_{j}^{f}(A)$ for the $j$-th Lipschitz-Killing
Curvatures of the set $A$ under the metric induced by the zero-mean Gaussian
field $f;$ for instance, in the case of spherical random eigenfunctions we
get immediately
\begin{equation*}
\mathcal{L}_{0}^{f_{\ell }}(\mathbb{S}^{2})=\mathcal{L}_{0}(\mathbb{S}_{%
\sqrt{\lambda _{\ell }/2}}^{2})=2\text{ , }\mathcal{L}_{1}^{f_{\ell }}(%
\mathbb{S}^{2})=0\text{ , }\mathcal{L}_{2}^{f_{\ell }}(\mathbb{S}^{2})=4\pi
\frac{\lambda _{\ell }}{2}\text{ .}
\end{equation*}%
For further notation, as in \cite{adlertaylor} we shall write
\begin{eqnarray*}
\rho _{j}(u) &:&=\frac{1}{(2\pi )^{1/2+j/2}}\exp (-u^{2}/2)H_{j-1}(u)\text{ , } j\geq 1 \\
\rho _{0}(u) &:&=1-\Phi (u)=\int_{u}^{\infty }\varphi(t)dt\text{ ,}
\end{eqnarray*}%
where as usual $\varphi(t)=(2\pi)^{-1/2}exp(-t^2/2)$ denotes the standard Gaussian density and we introduced the Hermite polynomials%
\begin{equation}
H_{k}(u):=(-1)^{k}\exp (\frac{u^{2}}{2})\frac{d^{k}}{du^{k}}\exp (-\frac{%
u^{2}}{2})\text{ , }k=0,1,2,...,u\in \mathbb{R}\text{ ;}  \label{Hermite}
\end{equation}%
for instance $H_{0}(u)=1,$ $H_{1}(u)=u,$ $H_{2}(u)=u^{2}-1,...$ Finally, we
shall introduce the \emph{flag coefficients}%
\begin{equation}
\left [ \begin{matrix} d \\ k \end{matrix} \right] :=\binom{d}{k}\frac{\omega _{d}}{\omega _{k}\omega _{d-k}}%
\text{ , }k=0,1,...,d\text{ .}  \label{flag}
\end{equation}

We are now in the position to state the following:

\begin{thm}
(Gaussian Kinematic Formula (\cite{Taylor}, \cite{adlertaylor}, Theorem
13.4.1)) Under Regularity Conditions, for all $j=0,1,...,n$ we have that%
\begin{equation}
\mathbb{E}\left[ \mathcal{L}_{j}^{f}(A_{u}(f;M))\right] =\sum_{k=0}^{d-j}
\left [ \begin{matrix} k+j \\ k \end{matrix} \right]
\rho _{k}(u)\mathcal{L}_{k+j}^{f}(M)\text{ .}  \label{GKF}
\end{equation}
\end{thm}

Before we proceed with some examples, it is worth discussing formula (\ref{GKF}). We are evaluating the expected value of a complex geometric
functional on a complicated excursion set, in very general circumstances
(under minimal regularity conditions on the field and on the manifold on
which it is defined). It is clear that the expected value should depend on
the manifold, on the threshold level, and on the field one considers, and
one may expect these three factors to be intertwined in a complicated
manner. On the contrary, formula (\ref{GKF}) shows that their role is
completely decoupled; more precisely

\begin{itemize}
\item the threshold $u$ enters the formula merely through the functions $%
\rho _{j}(u)$ which are very simple and fully universal (i.e., they do not
depend neither on the field nor on the manifold);

\item on the left-hand side Lipschitz-Killing Curvatures appear, but they
are computed on the original manifold, not on the excursion sets, and they
are therefore again extremely simple to compute;

\item the role of the field $f$ is confined to the new metric $g^{f}(.,.)$
that it induces and under which the Lipschitz-Killing Curvatures are
computed on both sides; under the (standard) assumption of isotropy, this
implies only a rescaling of the manifold by means of a factor depending only
on the derivative of the covariance function at the origin.
\end{itemize}

\begin{exa}
Let us consider a zero-mean isotropic Gaussian field $f$ defined on $\mathbb{%
S}^{d}$ (the unit sphere in $\mathbb{R}^{d+1});$ its covariance function can
be written as%
\begin{equation*}
\mathbb{E}\left[ f(x_{1})f(x_{2})\right] =\sum_{\ell =0}^{\infty }\frac{%
n_{\ell ,d}}{s_{d+1}}C_{\ell }G_{\ell;\frac{d}{2}}(\left\langle
x_{1},x_{2}\right\rangle )\text{ ,}
\end{equation*}%
where $s_{d+1}=(d+1)\omega_{d+1}$ is the surface measure of $\mathbb{%
S}^{d}$, $G_{\ell;\alpha}(.)$ denotes the normalized Gegenbauer polynomials of order $%
\alpha$, whereas
\begin{equation*}
n_{\ell ,d}=\frac{2\ell +d-1}{\ell }\binom{\ell +d-2}{\ell -1}\sim \frac{2}{%
(d-1)!}\ell ^{d-1},\text{ as }\ell \rightarrow \infty ,
\end{equation*}%
is the dimension of the eigenspace corresponding to the $\ell $-th
eigenvalue $\lambda _{\ell ;d}:=\ell (\ell +d-1);$ here $\left\{ C_{\ell
}\right\} $ is a sequence of non-negative weights which represent the
so-called angular power spectrum of the random field. The derivative of
the covariance function at the origin is%
\begin{equation*}
\mu :=\sum_{\ell =0}^{\infty }\frac{n_{\ell ,d}}{s_{d+1}}C_{\ell }\frac{%
\lambda _{\ell ;d}}{d}\text{ .}
\end{equation*}%
Recall the Lipschitz-Killing Curvatures of the manifold $\mathbb{S}_{\lambda
}^{d}:=\lambda \mathbb{S}^{d}$ are given by (\cite{adlertaylor}, page 179):%
\begin{equation*}
\mathcal{L}_{j}(\lambda \mathbb{S}^{d})=2\binom{d}{j}\frac{s_{d+1}}{s_{d+1-j}}%
\lambda ^{j}\text{ , }
\end{equation*}%
for $d-j$ even, and $0$ otherwise. Then the Gaussian Kinematic Formula
reads%
\begin{eqnarray*}
\mathbb{E}\left[ \mathcal{L}_{j}^{f}(A_{u}(f;\mathbb{S}^{d}))\right]
&=&\sum_{k=0}^{d-j}\rho _{k}(u)
\left [ \begin{matrix} k+j \\ k \end{matrix} \right]
\mathcal{L}_{k+j}(\sqrt{\mu}\mathbb{S}^{d}) \\
&=&\sum_{k=0}^{d-j}\rho _{k}(u)
\left [ \begin{matrix} k+j \\ k \end{matrix} \right]
\mathcal{L}_{k+j}(\mathbb{S%
}^{d})\mu^{(k+j)/2}.
\end{eqnarray*}
\end{exa}

\begin{exa}
As a special case of the previous example, assume $f=f_{\ell }$ is actually
a unit variance random eigenfunction on $\mathbb{S}^{2}$ corresponding to the eigenvalue $-\ell(\ell+1)$, $\ell=0,1,2,...$. Then the Gaussian
Kinematic Formula gives%
\begin{eqnarray*}
\mathbb{E}\left[ \mathcal{L}_{0}^{f_{\ell }}(A_{u}(f_{\ell };\mathbb{S}^{2}))%
\right] &=&\mathbb{E}\left[ \mathcal{L}_{0}(A_{u}(f_{\ell };\mathbb{S}^{2}))%
\right] \\
&=&2\left\{ 1-\Phi (u)\right\} +\frac{1}{2\pi }u\phi (u)\mathcal{(}4\mathcal{%
\pi )}\frac{\ell (\ell +1)}{2}\text{ ,}
\end{eqnarray*}%
\begin{equation*}
\mathbb{E}\left[ \mathcal{L}_{1}^{f_{\ell }}(A_{u}(f_{\ell };\mathbb{S}^{2}))%
\right] =\rho _{1}(u)
\left [ \begin{matrix} 2 \\ 1 \end{matrix} \right]
\mathcal{L}_{2}(\mathbb{S}^{2})\left\{ \frac{%
\ell (\ell +1)}{2}\right\}
\end{equation*}%
so that%
\begin{equation*}
\mathbb{E}\left[ \mathcal{L}_{1}(A_{u}(f_{\ell };\mathbb{S}^{2}))\right] = \pi \exp (-\frac{u^{2}}{2})\left\{ \frac{\ell (\ell +1)}{2}\right\}
^{1/2},
\end{equation*}%
and finally%
\begin{equation*}
\mathbb{E}\left[ \mathcal{L}_{2}(A_{u}(f_{\ell };\mathbb{S}^{2}))\right]
=\left\{ 1-\Phi (u)\right\} \mathcal{L}_{2}(\mathbb{S}^{2})=\left\{ 1-\Phi
(u)\right\} 4\pi \text{ }.
\end{equation*}
\end{exa}

\begin{exa}
In the special case of the nodal volume $\mathcal{L}_{d-1}(A_{0}(\mathbb{S}%
^{d}),f_{\ell })$ of random eigenfunctions, i.e., half the Hausdorff measure
of the zero-set of the eigenfunction, the Gaussian Kinematic Formula gives%
\begin{eqnarray*}
\mathbb{E}\left[ \mathcal{L}_{d-1}^{f}(A_{u}(f_{\ell };\mathbb{S}^{d}))%
\right] &=&(\frac{\lambda _{\ell }}{d})^{(d-1)/2}\mathbb{E}\left[ \mathcal{L}%
_{d-1}(A_{u}(f_{\ell };\mathbb{S}^{d}))\right] \\
&=&\rho _{1}(u)\frac{d\omega _{d}}{\omega _{1}\omega _{d-1}}\mathcal{L}_{d}(%
\mathbb{S}^{d})(\frac{\lambda _{\ell }}{d})^{d/2}
\end{eqnarray*}%
so that, recalling $\omega _{j}=\frac{\pi ^{j/2}}{\Gamma (\frac{j}{2}+1)}$
and $\mathcal{L}_{d}(\mathbb{S}^{d})=(d+1)\omega _{d+1}$%
\begin{eqnarray}
\mathbb{E}\left[ \mathcal{L}_{d-1}(A_{u}(f_{\ell };\mathbb{S}^{d}))\right]
&=&\frac{1}{2\pi }\exp (-\frac{u^{2}}{2})\frac{d\omega _{d}}{\omega
_{1}\omega _{d-1}}\mathcal{L}_{d}(\mathbb{S}^{d})(\frac{\lambda _{\ell }}{d}%
)^{1/2}  \notag \\
&=&\exp (-\frac{u^{2}}{2})\frac{\pi ^{d/2}}{\Gamma (\frac{d}{2})}(\frac{%
\lambda _{\ell }}{d})^{1/2}.  \label{ExpNodal}
\end{eqnarray}
\end{exa}

For $u=0$ (\ref{ExpNodal}) was derived for instance by \cite{Berard} (see
\cite{Wig}) and it is consistent with a celebrated conjecture by \cite{Yau},
which states that for $C^{\infty }$ manifolds the nodal volume of any eigenfunction
corresponding to the eigenvalue $E$ should belong to the interval $[c_{1}\sqrt{E},c_{2}\sqrt{E}]$
for some constants $0<c_{1}\leq c_{2}<\infty$.
The conjecture was settled for real analytic manifolds by \cite%
{DonnellyFefferman}; for smooth manifolds the lower bound was established
much more recently, see \cite{Logunov}, \cite{Logunov2}, \cite{Logunov3}
while the upper bound is addressed in \cite{Logunov4}. As a consequence of
the results in the next two Sections below in the case of the sphere in a
probabilistic sense the upper and lower constants can be taken nearly
coincident, in the limit of diverging eigenvalues.

\section{Wiener-Chaos Expansions, Variances and Correlations}

In view of the results detailed in the previous Section, the question
related to the expectation of intrinsic volumes in the case of Gaussian
fields can be considered completely settled. The next step of interest is
the computation of the corresponding variances, and the asymptotic laws of
fluctuations around the expected values, in the high-frequency regime. The
first rigorous results in this area can be traced back to a seminal paper by
Igor Wigman (\cite{Wig}) where the variance of the nodal length (i.e., $%
Len(f_{\ell },\mathbb{S}^{2}):=2\mathcal{L}_{1}(A_{0}(f_{\ell},\mathbb{S}^{2}))$ for random spherical harmonics in dimension 2 is computed and shown
to be asymptotic to
\begin{equation}
Var\left[ Len(f_{\ell },\mathbb{S}^{2})\right] =\frac{\log \ell }{32}%
+O_{\ell \rightarrow \infty }(1)\text{ }.  \label{Wig10}
\end{equation}%
We shall start instead from the derivation of variances and central limit
theorems for Lipschitz-Killing Curvatures of excursion sets at $u\neq 0,$
although these results were actually obtained more recently than (\ref{Wig10}).

Let us recall first the notion of Wiener chaos expansions. In the simplest
setting, consider $Y=G(Z)$ i.e., the transform of a zero mean, unit variance
Gaussian random variable $Z,$ such that $\mathbb{E}\left[ G(Z)^{2}\right]
<\infty ;$ it is well-known that the following expansion holds, in the $%
L^{2}(\Omega )$ sense:%
\begin{equation}
G(Z)=\sum_{q=0}^{\infty }\frac{J_{q}(G)}{q!}H_{q}(Z)\text{ ,}
\label{Chaos_exp}
\end{equation}%
where $\left\{ H_{q}(.)\right\} _{q=0,1,2,...}$ denotes the family of
Hermite polynomials that we introduced earlier in (\ref{Hermite}), and $%
J_{q}(G)$ are projection coefficients given by $J_{q}(G):=\mathbb{E}\left[
G(Z)H_{q}(Z)\right] $ (see i.e., \cite{Jansson}, \cite{noupebook})$.$ The
summands in (\ref{Chaos_exp}) are orthogonal, because when evaluated on pairs
of standard Gaussian variables $Z_{1},Z_{2},$ Hermite polynomials enjoy a
very simple formula for the computation of covariances:%
\begin{equation}
\mathbb{E}\left[ H_{q_{1}}(Z_{1})H_{q_{2}}(Z_{2})\right] =\delta
_{q_{1}}^{q_{2}}q_{1}!\left\{ \mathbb{E}\left[ Z_{1}Z_{2}\right]
\right\} ^{q_{1}},  \label{Wick}
\end{equation}%
where $\delta
_{q_{1}}^{q_{2}}$ denotes the Kronecker delta. Equation (\ref{Wick}) is just a special case of the celebrated \emph{Diagram
(or Wick's) Formula}, see \cite{noupebook} for much more discussion and
details. We thus have immediately%
\begin{equation*}
Var\left\{ G(Z)\right\} =\sum_{q=0}^{\infty }\frac{J_{q}^{2}(G)}{q!}\text{ .}
\end{equation*}%
More generally, let $\left\{ Z_{1},...,Z_{j},...\right\} $ be any array of
independent standard Gaussian variables, and consider elements of the form%
\begin{equation*}
H_{q_{1}}(Z_{1})\cdot ...\cdot H_{q_{p}}(Z_{p})\text{ , }q_{1}+...+q_{p}=q%
\text{ ;}
\end{equation*}%
the linear span (in the $L^{2}(\Omega )$ sense) of these random variables is
usually written $\mathcal{C}_{q}$ (denoted the $q$-th order Wiener chaos,
see again \cite{noupebook}) and we have the orthogonal decomposition%
\begin{equation*}
L^{2}(\Omega )=\bigoplus\limits_{q=0}^{\infty }\mathcal{C}_{q}\text{ .}
\end{equation*}

\subsection{Wiener-Chaos Expansions for Random Eigenfunctions}

Let us now explain how these techniques can be pivotal for the investigation
of fluctuations of geometric functionals. We start from the simplest case,
the excursion volume/area for the two-dimensional sphere, which we can write
as%
\begin{equation*}
\mathcal{L}_{2}(A_{u}(f_{\ell };\mathbb{S}^{2}))=\int_{\mathbb{S}^{2}}\mathbb{I}_{[u,\infty )}(f_{\ell
}(x))dx\text{ ,}
\end{equation*}%
$\mathbb{I}_{[u,\infty )}(.)$ denoting the indicator function of the
semi-interval $[u,\infty ).$ It is not difficult to show that%
\begin{eqnarray*}
J_{q}(\mathbb{I}_{[u,\infty )}(.)) &=&\mathbb{E}\left[ \mathbb{I}_{[u,\infty
)}(Z)H_{q}(Z)\right]  \\
&=&\int_{u}^{\infty }H_{q}(z)\phi (z)dz=(-1)^{q}H_{q-1}(u)\phi (u)\text{ ,}
\end{eqnarray*}%
the last result following by integration by parts, under the convention that
\begin{equation*}
(-1)H_{-1}(u)\phi (u):=1-\Phi (u)\text{ .}
\end{equation*}%
In view of (\ref{Chaos_exp}), we thus have (\cite{DI,MW})%
\begin{eqnarray*}
\mathcal{L}_{2}(A_{u}(f_{\ell };\mathbb{S}^{2})) &=&\int_{\mathbb{S}^{2}}\sum_{q=0}^{\infty
}(-1)^{q}H_{q-1}(u)\phi (u)\frac{H_{q}(f_{\ell }(x))}{q!}dx \\
&=&\sum_{q=0}^{\infty }\frac{(-1)^{q}}{q!}H_{q-1}(u)\phi (u)h_{\ell ;q}\text{
, where }h_{\ell ;q}=\int_{\mathbb{S}^{2}}H_{q}(f_{\ell }(x))dx\text{ ;}
\end{eqnarray*}%
as a consequence, we have also%
\begin{equation}
Var\left\{ \mathcal{L}_{2}(A_{u}(f_{\ell };\mathbb{S}^{2}))\right\} =\sum_{q=0}^{\infty }\frac{1}{%
(q!)^{2}}H_{q-1}^{2}(u)\phi ^{2}(u)Var\left\{ h_{\ell ;q}\right\} \text{ .}
\label{VarDefect}
\end{equation}

The crucial observation to be drawn at this stage is that the variances of
the components $\left\{ h_{\ell ;q}\right\} $ exhibit a form of phase
transition with respect to their order $q$, in the high-frequency/high
energy limit $\ell \rightarrow \infty $. In particular, a simple application
of the Diagram Formula (\ref{Wick}), isotropy and a change of variable yield%
\begin{eqnarray*}
Var\left\{ h_{\ell ;q}\right\} &=&\int_{\mathbb{S}^{2}\times \mathbb{S}^{2}}%
\mathbb{E}\left\{ H_{q}(f_{\ell }(x))H_{q}(f_{\ell }(y))\right\} dxdy \\
&=&8\pi ^{2}q!\int_{0}^{\pi }\left\{ P_{\ell }(\cos \theta )\right\}
^{q}\sin \theta d\theta \text{ ;}
\end{eqnarray*}%
for instance, for $q=2$ we obtain exactly%
\begin{equation*}
Var\left\{ h_{\ell ;q}\right\} =2\times 8\pi ^{2}\int_{0}^{\pi }P_{\ell
}^{2}(\cos \theta )\sin \theta d\theta =16\pi ^{2}\frac{2}{2\ell +1}\text{ .}
\end{equation*}%
Given two sequences of positive numbers ${a_n,b_n}$, we shall write $a_n \approx b_n$ when we have that $a_n/b_n \rightarrow c$ as $n \rightarrow \infty$, $c>0$. By means of the so-called Hilb's asymptotics (\cite{Szego},%
\cite{Wig}) it is possible to show that, as $\ell \rightarrow \infty $ (\cite%
{MW2014})%
\begin{eqnarray*}
Var\left\{ h_{\ell ;q}\right\} &\approx &\frac{1}{\ell ^{2}}\times
\int_{0}^{\ell \pi }\frac{1}{\psi ^{q/2}}\psi d\psi \\
&\approx &\left\{
\begin{array}{c}
\ell ^{-1}\text{ for }q=2 \\
\ell ^{-2}\log \ell\text{ for }q=4 \\
\ell ^{-2}\text{ for }q=3,5,...%
\end{array}%
\right. .
\end{eqnarray*}

Note that $h_{\ell ;1}\equiv 0$ for all $\ell =1,2,...,$ whereas the term
for $q=3$ requires an ad-hoc argument given in \cite{M2008, MW}. As a
consequence, the dominant terms in the variance expansion correspond to $q=2$
when $H_{1}(u)$ is non-zero, i.e., for $u\neq 0;$ for $u=0$ the even-order
chaoses vanish and all the remaining terms contribute by the same order of
magnitude with respect to $\ell $. In conclusion, we have that%
\begin{equation}
\mathcal{L}_{2}(A_{u}(f_{\ell };\mathbb{S}^{2}))-\mathbb{E}\left[ \mathcal{L}_{2}(A_{u}(f_{\ell };\mathbb{S}^{2}))\right] =\frac{%
1}{2}H_{1}(u)\phi (u)h_{\ell ;2}+O_{p}(\sqrt{\log \ell /\ell ^{2}})\text{ , }
\label{Area}
\end{equation}
and for $u\neq 0$%
\begin{equation*}
Var\left\{ \mathcal{L}_{2}(A_{u}(f_{\ell };\mathbb{S}^{2}))\right\} \sim \left\{ \frac{1}{2}%
H_{1}(u)\phi (u)\right\} ^{2}Var\left\{ h_{\ell ;2}\right\} \text{ , as }%
\ell \rightarrow \infty \text{ .}
\end{equation*}%
Because
\begin{equation*}
h_{\ell ;2}=\int_{\mathbb{S}^{2}}\left\{ f_{\ell }^{2}(x)-1\right\}
dx=\left\Vert f_{\ell }\right\Vert _{L^{2}(\mathbb{S}^{2})}^{2}-\mathbb{E}%
\left[ \left\Vert f_{\ell }\right\Vert _{L^{2}(\mathbb{S}^{2})}^{2}\right]
\text{ ,}
\end{equation*}%
equation (\ref{Area}) is basically stating that the fluctuations in the
excursion area for $u\neq 0$ are dominated by the fluctuations in the random
norm of the eigenfunctions.

Interestingly, the same behaviour characterizes also the other
Lipschitz-Killing Curvatures; for the boundary length we have the expansion%
\begin{equation*}
2\mathcal{L}_{1}(A_{u}(f_{\ell };\mathbb{S}^{2}))=\lim_{\varepsilon
\rightarrow 0}\int_{\mathbb{S}^{2}}\left\Vert \nabla f_{\ell }(x)\right\Vert
\delta _{\varepsilon }(f_{\ell }(x)-u)dx
\end{equation*}%
which holds both $\omega $-almost surely and in $L^{2}(\Omega )$; here we write $\delta _{\varepsilon }(\cdot)=\frac{1}{2\varepsilon} \mathbb{I}(\cdot)$.
Similarly for the Euler-Poincar\'{e} Characteristic we have
\begin{equation*}
\mathcal{L}_{0}(A_{u}(f_{\ell };\mathbb{S}^{2}))=\lim_{\varepsilon
\rightarrow 0}\int_{\mathbb{S}^{2}}\det \left\{ \nabla ^{2}f_{\ell
}(x)\right\} \delta _{\varepsilon }(\nabla f_{\ell }(x))\mathbb{I}%
_{[u,\infty )}\mathbb{(}f_{\ell }(x))dx\text{ .}
\end{equation*}

Similar arguments can be developed, expanding the integrand function into
polynomials evaluated on the random vectors $\left\{ \nabla ^{2}f_{\ell
}(.),\nabla f_{\ell }(.),f_{\ell }(.)\right\}$; algebraic simplifications
occur and the expansions read as follows:

\begin{thm}
\label{2GKF} As $\ell \rightarrow \infty $, for $j=0,1,2$%
\begin{eqnarray}
&&\mathcal{L}_{j}(A_{u}(f_{\ell },\mathbb{S}^{2}))-\mathbb{E}\left[ \mathcal{%
L}_{j}(A_{u}(f_{\ell };\mathbb{S}^{2}))\right]   \label{Leading} \\
&=&-\frac{1}{2}
\left [ \begin{matrix} 2 \\ 2-j \end{matrix} \right]
u\rho _{2-j}^{\prime }(u)\left( \lambda
_{\ell }/2\right) ^{(2-j)/2}\int_{\mathbb{S}^{2}}H_{2}(f_{\ell
}(x))dx+R_{\ell ;j}\text{ ,}  \notag
\end{eqnarray}%
where
\[ 
\mathbb{E}[R^2_{\ell ;j}]=o_{\ell
\rightarrow \infty }(\ell^{3-2j})\text{ ;}
\]
as a consequence, we have also the following Variance asymptotics%
\[
Var\left\{ \mathcal{L}_{j}(A_{u}(f_{\ell };\mathbb{S}^{2}))\right\}
\]
\begin{equation}
=\frac{1}{4}\left\{\left [ \begin{matrix} 2 \\ 2-j \end{matrix} \right]
u\rho _{2-j}^{\prime }(u)\left( \lambda _{\ell
}/2\right) ^{(2-j)/2}\right\} ^{2}\times \frac{32\pi ^{2}}{2\ell +1}+o_{\ell
\rightarrow \infty }(\lambda _{\ell }^{2-j-1})\text{ .}  \label{Variance}
\end{equation}
\end{thm}

Some features of the previous result are worth discussing:

\begin{itemize}
\item The asymptotic behaviour of all the Lipschitz-Killing Curvatures is
proportional to a sequence of scalar random variables $\left\{ h_{\ell
;2}\right\} _{\ell \in \mathbb{N}}.$ As a consequence, these geometric
functionals are fully correlated in the high-energy limit $\ell \rightarrow
\infty ;$

\item For the same reasons, these functionals are also fully correlated, in
the high energy limit, when evaluated across different levels $u_{1},u_{2}$: 
for the boundary length, this correlation phenomenon was first noted by \cite{WigSur};

\item The leading terms all disappear in the "nodal" case $u=0$ where the
variances are hence an order of magnitude smaller. This is an instance of
the so-called Berry's cancellation phenomenon (\cite{Wig}), to which we
shall return in the following Section. We noted before that the leading
terms are proportional to the centred random norm; it is thus natural that
these terms should disappear in the nodal case, which is independent of
scaling factors. Note that for $j=0$ the cancellation of the leading term occurs also at $u=1$.

\end{itemize}

\begin{rem}
The proof of Theorem \ref{2GKF} was given in \cite{CM2018}, in the case of
the $2$-dimensional sphere $\mathbb{S}^{2}.$ However, we conjecture the
result to hold as stated for spherical eigenfunctions in arbitrary
dimension, see below for more details. Extensions have also been given to
cover for instance the two-dimensional torus (see \cite{CMR}), for which a
formula completely analogous to (\ref{2GKF}) holds.
\end{rem}

Similar results can be shown to hold for other geometric functionals; let us
consider for instance critical values, defined by%
\begin{equation*}
\mathcal{N}_{u}(f_{\ell };\mathbb{S}^{2})=\#\left\{ x\in \mathbb{S}%
^{2}:\nabla f_{\ell }(x)=0\text{ and }f_{\ell }(x)\geq u\right\} \text{ .}
\end{equation*}%
The asymptotic variance of $\left\{ \mathcal{N}_{u}(f_{\ell };\mathbb{S}%
^{2})\right\} _{\ell =1,2,...}$ was established in \cite{CMW}, \cite{CW},
and in particular we have%
\begin{eqnarray*}
\mathbb{E}\left[ \mathcal{N}_{u}(f_{\ell };\mathbb{S}^{2})\right]
&=&\lambda _{\ell }g_{1}(u)\text{ , } \\
g_{1}(u) &=&\frac{1}{\sqrt{2\pi }}\int_{u}^{\infty
}(2e^{-t^{2}}+(t^{2}-1)e^{-t^{2}/2})dt \\
&=&u\phi (u)+\sqrt{2}(1-\Phi (\sqrt{2}u))\text{ ,}
\end{eqnarray*}%
\begin{eqnarray*}
Var\left[ \mathcal{N}_{u}(f_{\ell };\mathbb{S}^{2})\right]  &=&\frac{1}{4}%
\lambda _{\ell }^{2}g_{2}^{2}(u)Var\left\{ \int_{\mathbb{S}%
^{2}}H_{2}(f_{\ell }(x))dx\right\} +o_{\ell \rightarrow \infty }(\ell ^{3})
\\
&=&\frac{1}{4}\lambda _{\ell }^{2}g_{2}^{2}(u)\frac{2(4\pi )^{2}}{2\ell +1}%
+o_{\ell \rightarrow \infty }(\ell ^{3})\text{ ,}
\end{eqnarray*}
where
\begin{equation*}
g_{2}(u)=\int_{u}^{\infty }\frac{1}{\sqrt{8\pi }}%
e^{-3t^{2}/2}(2-6t^{2}-e^{-t^{2}}(1-4t+t^{4}))dt\text{ .}
\end{equation*}%
Later in \cite{CM2020} it was shown that the critical values above the
threshold level $u$ satisfy the asymptotic%
\begin{eqnarray*}
&&\mathcal{N}_{u}(f_{\ell };\mathbb{S}^{2})-\mathbb{E}\left[ \mathcal{N}%
_{u}(f_{\ell };\mathbb{S}^{2})\right]  \\
&=&\frac{1}{2}\lambda _{\ell }g_{2}(u)\int_{\mathbb{S}^{2}}H_{2}(f_{\ell
}(x))dx+o_{p}(\sqrt{Var\left[ \mathcal{N}_{u}(f_{\ell };\mathbb{S}^{2})%
\right] )}\text{ ,}
\end{eqnarray*}%
As a consequence, one has also, for all $u\neq 0,1$ the following
correlation result%
\begin{eqnarray*}
&&Corr^{2}\left\{ \mathcal{N}_{u}(f_{\ell };\mathbb{S}^{2}),\mathcal{L}%
_{j}(A_{u}(f_{\ell };\mathbb{S}^{2}))\right\}  \\
&:&=\frac{Cov^{2}\left\{ \mathcal{N}_{u}(f_{\ell };\mathbb{S}^{2}),\mathcal{L%
}_{j}(A_{u}(f_{\ell };\mathbb{S}^{2}))\right\} }{Var\left\{ \mathcal{N}%
_{u}(f_{\ell };\mathbb{S}^{2})\right\} Var\left\{ \mathcal{L}_{j}(A_{u}(f_{\ell };\mathbb{S}^{2}))\right\} }\rightarrow 1\text{ , as }\ell
\rightarrow \infty \text{ ;}
\end{eqnarray*}%
the value $u=1$ has to be excluded only for $j=0$. We also have that
\begin{equation*}
Corr^{2}\left\{ \mathcal{N}_{u_{1}}(f_{\ell };\mathbb{S}^{2}),\mathcal{N}%
_{u_{2}}(f_{\ell };\mathbb{S}^{2})\right\} \rightarrow 1\text{ , as }\ell
\rightarrow \infty \text{ ,}
\end{equation*}%
that is, asymptotically full correlation between the number of critical
values above any two non-zero thresholds $u_{1},u_{2}.$

As for the Lipschitz-Killing Curvatures, a form of Berry's cancellation
occurs at $u=0$ and $u \rightarrow \pm \infty ;$ the total number of critical points has
then a lower-order variance (see \cite{CW}), as we shall discuss in the next
section.

\subsection{Quantitative Central Limit Theorems}

The results reviewed in the previous subsection can be considered as
following from a \emph{Reduction Principle }(see \cite{DehTaq}), where the
limiting behaviour of $\left\{ \mathcal{N}_{u}(f_{\ell };\mathbb{S}^{2}),%
\mathcal{L}_{j}(A_{u}(f_{\ell };\mathbb{S}^{2})\right\} $ is dominated by a
deterministic function of the threshold level $u,$ times a sequence of
random variables $\left\{ h_{\ell ;2}\right\} $ which do not depend on $u.$
To derive the asymptotic law of these fluctuations, it is hence enough to
investigate the convergence in distribution of $\left\{ h_{\ell ;2}\right\} ,
$ as $\ell \rightarrow \infty $. In fact, it is possible to show a stronger
result, namely a \emph{Quantitative Central Limit Theorem}; to this aim, let
us recall that the Wasserstein distance between two random variables $X$ and
$Y$ is defined by%
\begin{equation*}
d_{W}(X,Y):=\sup_{h\in Lip(1)}\left\vert \mathbb{E}h(X)-\mathbb{E}%
h(Y)\right\vert \text{ ,}
\end{equation*}%
where $Lip(1)$ denotes the class of Lipschitz functions of constant 1, i.e.,
$\left\vert h(x)-h(y)\right\vert \leq \left\vert x-y\right\vert $ for all $%
x,y\in \mathbb{R}$. $D_{W}(.,.)$ defines a metric on the space of
probability distributions (for more details and other examples of
probability metrics, see \cite{noupebook}, Appendix C). Taking $Z\sim N(0,1)$
to be a standard Gaussian random variable, a Quantitative Central Limit
theorem is defined as a result of the form%
\begin{equation*}
\lim_{n\rightarrow \infty }d_{W}(\frac{X_{n}-\mathbb{E}X_{n}}{\sqrt{%
Var(X_{n})}},Z)=0\text{ .}
\end{equation*}%
The field of Quantitative Central Limit Theorems has been very active in the
last few decades; more recently, a breakthrough has been provided by the
discovery of the so-called \emph{Stein-Malliavin approach} by
Nourdin-Peccati (\cite{nuape, NP09, noupebook}). These results entail that
for sequences of random variables belonging to a Wiener-chaos, say $\mathcal{C}_q$, a
quantitative central limit theorem for the Wasserstein distance can be given
simply controlling the fourth-moment of $X_{n},$ as follows:%
\begin{equation}
d_{W}(\frac{X_{n}-\mathbb{E}X_{n}}{\sqrt{Var(X_{n})}},Z)\leq \sqrt{\frac{2q-2%
}{3\pi q}}\sqrt{\mathbb{E}\left[ (\frac{X_{n}-\mathbb{E}X_{n}}{\sqrt{%
Var(X_{n})}})^{4}\right] -3}\text{ .}  \label{FMB}
\end{equation}%
Similar results hold for other probability metrics, for instance the
Kolmogorov and Total Variation distances, see again \cite{noupebook}.

Quantitative Central Limit Theorems lend themselves to an immediate
application for the sequences $\left\{ h_{\ell ;q}\right\} $ that we
introduced above. It should be noted indeed that by construction all these
random variables belong to the $q$-th order Wiener chaos; it is then
possible to exploit (\ref{FMB}) to obtain Quantitative Central Limit Theorems
for these polyspectra at arbitrary orders: their fourth moment can be
computed by means of the Diagram formula. These results were first given in
\cite{MW} and then refined in \cite{MR2015}, yielding the following

\begin{thm}
\label{QCLT_h} As $\ell \rightarrow \infty $%
\begin{equation*}
d_{W} \left (\frac{h_{\ell ;q}-\mathbb{E}\left[ h_{\ell ;q}\right] }{\sqrt{Var(h_{\ell ;q})}}, Z \right )=\left\{
\begin{array}{c}
O(\frac{1}{\sqrt{\ell }})\text{ for }q=2,3 \\
O(\frac{1}{\log \ell })\text{ for }q=4 \\
O(\ell ^{-1/4})\text{ for }q=5,6,...%
\end{array}%
\right. .
\end{equation*}
\end{thm}

Now, we have just shown that for nonzero thresholds $u\neq 0$ the
Lipschitz-Killing Curvatures and the critical values are indeed proportional
to a term belonging to the second-order chaos, plus a remainder that it is
asymptotically negligible. The following Quantitative Central Limit Theorem
then follows immediately (see \cite{MW}, \cite{RossiJTP}, \cite{CM2018}).

\begin{thm}
As $\ell \rightarrow \infty ,$ for $u\neq 0$ ($j=1,2)$ and for $u\neq 0,1$
(for $j=0$) we have that%
\begin{equation*}
d_{W}(\frac{\mathcal{L}_{j}(A_{u}(f_{\ell };\mathbb{S}^{2}))-\mathbb{E}\left[
\mathcal{L}_{j}(A_{u}(f_{\ell };\mathbb{S}^{2}))\right] }{\sqrt{Var(\mathcal{%
L}_{j}(A_{u}(f_{\ell };\mathbb{S}^{2}))}},Z)=O(\ell ^{-1/2})\text{ .}
\end{equation*}
\end{thm}

\subsection{A Higher-Dimensional Conjecture}

The results we discussed so far have been limited to random-spherical
harmonics on the two-dimensional sphere $\mathbb{S}^{2}.$ Research in
progress suggests however that further generalizations should hold: to this
aim, let us define the set of singular points $P_{j}:=\left\{ u\in \mathbb{R}%
:u\rho _{j}^{\prime }(u)=0\right\} $ (for instance, $P_{0}=P_{1}=\left\{
0\right\} ,$ $P_{2}=\left\{ 0,1\right\} ,$ $P_{3}=\left\{ 0,\pm \sqrt{3}%
\right\} ,...).$ Let us now consider Gaussian random eigenfunctions on the
higher-dimensional unit sphere $\mathbb{S}^{d}$, e.g.%
\begin{equation*}
\Delta _{\mathbb{S}^{d}}f_{\ell ;d}=-\lambda _{\ell ;d}f_{\ell ;d}\text{ , }%
\lambda _{\ell ;d}:=\ell (\ell +d-1)\text{ ;}
\end{equation*}%
these eigenfunctions are normalized so that (see \cite{MR2015},\cite%
{Rossi2018})%
\begin{equation*}
\mathbb{E}\left[ f_{\ell ;d}\right] =0\text{ , }\mathbb{E}\left[ f_{\ell
;d}^{2}\right] =1\text{ , }\mathbb{E}\left[ f_{\ell ;d}(x)f_{\ell ;d}(y)%
\right] =G_{\ell ;d/2}(\left\langle x,y\right\rangle )\text{ ,}
\end{equation*}%
where as before $G_{\ell ;d/2}(.)$ is the standardized $\ell $-th Gegenbauer polynomial
of order $\frac{d}{2}$ (normalized with $G_{\ell ;d/2}(1)=1$); it is convenient to recall that%
\begin{equation*}
G_{\ell ;d/2}^{\prime }(1)=\frac{\lambda _{\ell ;d}}{d}\text{ .}
\end{equation*}%
We recall also that the dimension of the corresponding eigenspaces is
\begin{equation*}
n_{\ell ;d}=\frac{2\ell +d-1}{\ell }\binom{\ell +d-2}{\ell -1}\sim \frac{2}{%
(d-1)!}\ell ^{d-1},\text{ as }\ell \rightarrow \infty \text{ .}
\end{equation*}%
By means of Parseval's equality we have also as a consequence%
\begin{eqnarray*}
Var\left[ \int_{\mathbb{S}^{d}}H_{2}(f_{\ell ;d}(x))dx\right]  &=&\frac{%
2s_{d}^{2}}{n_{\ell ;d}}=\frac{2(d+1)^{2}\omega _{d+1}^{2}}{n_{\ell ;d}} \\
&\sim &\frac{(d+1)^{2}\omega _{d+1}^{2}(d-1)!}{\ell ^{d-1}}\text{ as }\ell
\rightarrow \infty \text{ .}
\end{eqnarray*}%
We then propose the following

\begin{conjecture}
\label{H2GKF} As $\ell \rightarrow \infty ,$ for all $k=0,1,...,d$ we have
that%
\begin{equation*}
\mathcal{L}_{k}(A_{u}(f_{\ell };\mathbb{S}^{d}))-\mathbb{E}\left[ \mathcal{L}%
_{k}(A_{u}(f_{\ell };\mathbb{S}^{d}))\right]
\end{equation*}%
\begin{equation*}
=-\frac{1}{2}\left [ \begin{matrix} d \\ k \end{matrix} \right]
\rho _{d-k}^{\prime }(u)u\left( \frac{\lambda
_{\ell ;d}}{d}\right) ^{(d-k)/2}\int_{\mathbb{S}^{d}}H_{2}(f_{\ell
;d}(x))dx+o_{p}(\sqrt{\ell ^{d-2k+1}})\text{ .}
\end{equation*}
\end{conjecture}

\begin{rem}
An immediate consequence of this conjecture would be%
\begin{eqnarray*}
\frac{\mathcal{L}_{k}(A_{u}(f_{\ell };\mathbb{S}^{d}))-\mathbb{E}\left[
\mathcal{L}_{k}(A_{u}(f_{\ell };\mathbb{S}^{d}))\right] }{\sqrt{Var\left[
\mathcal{L}_{k}(A_{u}(f_{\ell };\mathbb{S}^{d}))\right] }} &=&\frac{h_{\ell
;q}}{\sqrt{Var\left[ h_{\ell ;d}(2)\right] }}+o_{p}(1)\text{ ,} \\
h_{\ell ;q} &=&\int_{\mathbb{S}^{d}}H_{2}(f_{\ell ;d}(x))dx\text{ .}
\end{eqnarray*}
\end{rem}

\begin{rem}
The remainder term in Conjecture (\ref{H2GKF}) is expected to be $O(\sqrt{\ell ^{d-2k}}),$
in the $L^{2}(\Omega )$ sense.
\end{rem}

Three further consequences of Conjecture \ref{H2GKF} would be the following:

\begin{itemize}
\item \emph{(Variance Asymptotics)} As $\ell \rightarrow \infty ,$ for all $k=0,1,...,d$ and for non-singular points $u \notin P_{d-k}$,%
\begin{equation*}
Var\left\{ \mathcal{L}_{k}(A_{u}(f_{\ell };\mathbb{S}^{d}))\right\} 
\end{equation*}%
\begin{equation*}
=\frac{H_{d-k}^{2}(u)\phi ^{2}(u)u^{2}}{(2\pi d)^{(d-k)}}\frac{d!}{(d-k)!k!}%
\frac{\omega _{d}^{2}\omega _{d+1}^{2}}{\omega _{k}^{2}\omega _{d-k}^{2}}%
\frac{(d+1)^{2}\lambda _{\ell ;d}^{d-k}}{2n_{\ell ;d}}+o(\ell ^{d-2k+1})%
\text{ .}
\end{equation*}

\item \emph{(Central Limit Theorem) }As $\ell \rightarrow \infty ,$ for all $k=0,1,...,d$ and for non-singular points $u \notin P_{d-k}$,
\begin{equation*}
d_{W}\left( \frac{\mathcal{L}_{k}(A_{u}(f_{\ell };\mathbb{S}^{d}))-\mathbb{E}%
\left[ \mathcal{L}_{k}(A_{u}(f_{\ell };\mathbb{S}^{d}))\right] }{\sqrt{Var%
\left[ \mathcal{L}_{k}(A_{u}(f_{\ell };\mathbb{S}^{d}))\right] }},Z\right) =o(1)\text{ ,}
\end{equation*}
where $Z \sim \mathcal{N}(0,1)$.
\item \emph{(Correlation Asymptotics)} As $\ell \rightarrow \infty ,$ for
all $k_{1},k_{2}=0,1,...,d$ and all $u_{1},u_{2}$ such that $%
u_{1}u_{2}H_{d-k_{1}}(u_{1})H_{d-k_{2}}(u_{2})\neq 0$%
\begin{equation*}
\lim_{\ell \rightarrow \infty }Corr^{2}\left( \mathcal{L}_{k_{1}}(A_{u}(f_{\ell };\mathbb{S}^{d})),\mathcal{L}_{k_{2}}(A_{u}(f_{\ell };\mathbb{S}^{d}))\right) =1\text{ .}
\end{equation*}
\end{itemize}

The driving rationale behind these conjectures is the \emph{ansatz} that the
asymptotic variance of the geometric functionals should be governed by
fluctuations in the random $L^{2}(\mathbb{S}^{d})$ norm of the
eigenfunctions, for non-singular points $u\notin P_{j}$. In this sense, we
believe the result has even greater applicability, for instance to cover
combinations of random eigenfunctions defined on more general
submanifolds of $\mathbb{R}^{n},$ such as \emph{Berry Random Waves} or
"short windows" averages of isotropic random eigenfunctions on general
manifolds (see \cite{Berry 1977}, \cite{Berry 2002}, \cite{Zel}, \cite{npr}, \cite{Duerickx}, \cite{DalmaoEstradeLeon}). These issues are the object of
currently ongoing research.

\section{Nodal Cases: Berry Cancellation and the Role of The Fourth-Order
Chaos}

The previous Section has discussed the behaviour of geometric functionals
for non-zero threshold levels $u\neq 0;$ under isotropy, it has been shown
that all these functionals are asymptotically proportional, in the $%
L^{2}(\Omega )$ sense, to a single random variable representing the
(centred) random $L^{2}(\mathbb{S}^{2})$-norm of the eigenfunction. This
dominant term has been shown to disappear in the nodal case $u=0$ (and, more
generally, for $\rho _{d-k}^{\prime }(u)u=0,$ i.e., for the singular points $%
u\in P_{j}$); the asymptotic behaviour must then be derived by a different
route in these circumstances.

As mentioned above, the first paper to investigate the variance of the nodal
length for random spherical harmonics was the seminal work by Igor Wigman (%
\cite{Wig}), which made rigorous an \emph{ansatz} by Michael Berry in the
Physical literature (\cite{Berry 2002}). In particular, by using an
higher-order version of the Expectation Metatheorem (see again \cite%
{adlertaylor}, \cite{azaiswschebor}) the following representation for the
second moment of the nodal length can be given:%
\begin{equation*}
\mathbb{E}\left[ \left\{ Len(f_{\ell };\mathbb{S}^{2})\right\}^2 \right]
\end{equation*}%
\begin{equation*}
=\int_{\mathbb{S}^{2}\times \mathbb{S}^{2}}\mathbb{E}\left[ \left.
\left\Vert \left\{ \nabla f_{\ell }(t_{1})\right\} \right\Vert \left\Vert
\left\{ \nabla f_{\ell }(t_{2})\right\} \right\Vert \right\vert f_{\ell
}(t_{1})=0,f_{\ell }(t_{2})=0\right]
\end{equation*}%
\[
\times p_{f_{\ell }(t_{1}),f_{\ell
}(t_{2})}(0,0)\sigma _{g}(dt_{1})\sigma _{g}(dt_{2})\text{ },
\]
where as before we write $Len(f_{\ell };\mathbb{S}^{2})=2\mathcal{L}_{1}(A_{0}(f_{\ell };\mathbb{S}^{2}))$ for the nodal length. The integrand in the previous formula is denoted the \emph{2-point
correlation function of the nodal length }and generalizes the Kac-Rice
argument to second-order moments; analogous generalizations are possible for
the other geometric functionals we considered and for higher-order moments
as well (see \cite{adlertaylor}). By means of a challenging and careful
expansion of this correlation function and a deep investigation of its
behaviour for $\ell \rightarrow \infty ,$ Wigman was able to investigate the
asymptotic for the variance of the nodal length and to show that (\ref{Wig10})
holds.

A natural question which was investigated shortly after this seminal paper
was the possibility to derive the asymptotic variances of nodal statistics,
and further characterizations such as the law of the asymptotic
fluctuations, in terms of the Wiener-Chaos expansions that we discussed in
the previous Section. The first efforts were devoted to the analysis of the
"nodal area" $\mathcal{L}_{2}(A_{0}(f_{\ell };\mathbb{S}^{2})),$ for which
it is easily shown that all even-order terms vanish at $u=0;$ from (\ref%
{VarDefect}) we are then left with (see \cite{MW2014})%
\begin{equation*}
Var\left\{ \mathcal{L}_{2}(A_{0}(f_{\ell };\mathbb{S}^{2}))\right\} =\frac{1%
}{\ell ^{2}}\sum_{q=1}^{\infty }\frac{c_{2q+1}}{2\pi q!}H_{2q}^{2}(0)+o(\ell
^{-2})\text{ ,}
\end{equation*}%
where%
\begin{eqnarray*}
c_{2q+1} &=&\lim_{\ell \rightarrow \infty }\ell ^{2}\int_{0}^{\pi }P_{\ell
}^{2q+1}(\cos \theta )\sin \theta d\theta  \\
&=&\int_{0}^{\infty }J_{0}^{2q+1}(\psi )\psi d\psi \text{ , }J_{0}(\psi
):=\sum_{k=0}^{\infty }\frac{(-1)^{k+1}(x/2)^{2k}}{(k!)^{2}}\text{ .}
\end{eqnarray*}%
The computation of the variance and the results in Theorem \ref{QCLT_h}
lead easily also to a Central Limit Theorem, which was given first in \cite%
{MW} and then extended to higher dimensions by \cite{RossiJTP}.

\begin{thm}
(\cite{MW}) As $\ell \rightarrow \infty $%
\begin{equation*}
d_{W}(\frac{\mathcal{L}_{2}(A_{0}(f_{\ell };\mathbb{S}^{2}))-\mathbb{E}\left[
\mathcal{L}_{2}(A_{0}(f_{\ell };\mathbb{S}^{2}))\right] }{\sqrt{Var\left\{
\mathcal{L}_{2}(A_{0}(f_{\ell };\mathbb{S}^{2}))\right\} }},Z)=o(1)\text{ ,}
\end{equation*}%
and hence
\begin{equation*}
\frac{\mathcal{L}_{2}(A_{0}(f_{\ell };\mathbb{S}^{2}))-\mathbb{E}\left[
\mathcal{L}_{2}(A_{0}(f_{\ell };\mathbb{S}^{2})\right] }{\sqrt{Var\left\{
\mathcal{L}_{2}(A_{0}(f_{\ell };\mathbb{S}^{2}))\right\} }}\rightarrow
_{d}\mathcal{N}(0,1))\text{ .}
\end{equation*}
\end{thm}

The proof of the previous result is standard; in short, the idea is to write%
\begin{equation*}
\mathcal{L}_{2}(A_{0}(f_{\ell };\mathbb{S}^{2}))-\mathbb{E}\left[ \mathcal{L}%
_{2}(A_{0}(f_{\ell };\mathbb{S}^{2})\right] =\sum_{k=1}^{M}\frac{(-1)^{2k+1}%
}{(2k+1)!}H_{2k}(u)\phi (u)h_{\ell ;2k+1}+R_{M}\text{ ,}
\end{equation*}%
where the remainder term is such that, as $M\rightarrow \infty ,$%
\begin{equation*}
R_{M}=\sum_{k=M+1}^{\infty }\frac{(-1)^{2k+1}}{(2k+1)!}H_{2k}(u)\phi
(u)h_{\ell ;2k+1}=o_{p}(\sqrt{Var\left\{ \mathcal{L}_{2}(A_{0}(f_{\ell };\mathbb{S}^{2}))\right\} })\text{ .}
\end{equation*}%
It is then enough to show that the Central Limit Theorem holds for $M$
(sufficiently large but) finite; this can be achieved by an application of
the multivariate Fourth Moment Theorem to the terms $(h_{\ell
;3},...,h_{\ell ;2M+1})$ (see \cite{noupebook})$.$ It should be noted that
in the case of the Defect the limiting behaviour depends on the full
sequence $\left\{ h_{\ell ;2k+1}\right\} _{k=1,2,...};$ this is due to the
exact disappearance of the two natural candidates to be leading terms, that
is, $\left\{ h_{\ell ;2}\right\} $ and $\left\{ h_{\ell ;4}\right\} ,$ both
whose coefficients vanish for $u=0.$

It is thus even more remarkable that for the nodal lines the situation
simplifies drastically, to yield the following result.

\begin{thm}
(\cite{MRW}) \label{nodal} As $\ell \rightarrow \infty $%
\begin{equation}
Len(f_{\ell };\mathbb{S}^{2})-\mathbb{E}\left[ Len(f_{\ell };\mathbb{S}^{2})%
\right] =-\frac{1}{4}\sqrt{\frac{\lambda _{\ell }}{2}}\frac{1}{4!}h_{\ell
;4}+o_{p}(\sqrt{Var\left\{ h_{\ell ;4}\right\} })\text{ ,}  \label{MRW}
\end{equation}
and hence, in view of (\ref{QCLT_h})%
\begin{equation*}
d_{W}(\frac{Len(f_{\ell };\mathbb{S}^{2})-\mathbb{E}\left[ Len(f_{\ell };\mathbb{S}^{2})\right] }{\sqrt{Var\left\{ Len(f_{\ell};\mathbb{S}^{2})\right\} }},Z)=o(1)\text{ .}
\end{equation*}
\end{thm}

The most notable aspect of Theorem \ref{nodal} is that the limiting
behaviour of nodal lines is asymptotically fully correlated with the
sequence of random variables $\left\{ h_{\ell ;4}\right\}$, so that in
principle it would be possible to "predict" nodal lengths by simply
computing the integral of a fourth-order polynomial of the eigenfunctions
over the sphere.

A natural question that arises is the structure of correlation among
functionals evaluated at different thresholds and those considered for the
nodal case $u=0$. Focussing for instance on the boundary length, it is
immediate to understand that the latter, which is dominated by the second
order chaos term $\left\{ h_{\ell ;2\text{ }}\right\}$ when $u\neq 0,$ must
be independent from the nodal length, which is asymptotically proportional
to $\left\{ h_{\ell ;4}\right\} .$ A more refined analysis, however, should
take into account the fluctuations of the boundary length when the effects
of the random norm $\left\Vert f_{\ell }\right\Vert _{L^{2}(\mathbb{S}^{2})}$
is subtracted, that is, dropping the second-order chaos term from the Wiener expansion. This corresponds to the evaluation of the so-called partial
correlation coefficients $Corr^{\ast }$, for which it was shown in \cite%
{MR21} that
\begin{equation*}
\lim_{\ell \rightarrow \infty }Corr^{\ast }(Len(f_{\ell };\mathbb{S}^{2}),%
\mathcal{L}_{1}(A_{u}(f_{\ell };\mathbb{S}^{2})))=1\text{ .}
\end{equation*}%
More explicitly, when compensating the effect of random norm fluctuations
the boundary length at any threshold $u\neq 0$ can be fully predicted on the
basis of the knowledge of the nodal length, up to a remainder term which is
asymptotically negligible in the limit $\ell \rightarrow \infty .$ It is
interesting to note that a similar phenomenon occurs also for the total
number of critical points, for which (building on earlier computations by
\cite{CW}) it was shown in (\cite{CM2021}) that

\begin{equation*}
\mathcal{N}_{-\infty }(f_{\ell };\mathbb{S}^{2})-\mathbb{E}\left[ \mathcal{N}%
_{-\infty }(f_{\ell };\mathbb{S}^{2})\right] =-\frac{\lambda _{\ell }}{%
2^{3}3^{2}\sqrt{3}\pi }h_{\ell ;4}+o_{p}(\ell ^{2}\log \ell )\ ;
\end{equation*}%
as a consequence, the nodal length of random spherical harmonics and the
number of their critical points are perfectly correlated in the high-energy
limit:%
\begin{equation*}
\lim_{\ell \rightarrow \infty }Corr^{2}(Len(f_{\ell };\mathbb{S}^{2}),%
\mathcal{N}_{-\infty }(f_{\ell };\mathbb{S}^{2}))=1\text{ .}
\end{equation*}
Let us now denote by $Len^{\ast }(u)$ the boundary length at level $u$ after
the fluctuations induced by the random norm have been subtracted (e.g.,
after removing its projection on the second-order chaos); moreover, for
brevity's sake we write%
\begin{eqnarray*}
\mathcal{L}_{j}(A_{u}(f_{\ell };\mathbb{S}^{2})) &=&\mathcal{L}_{j}(u)\text{
, }j=0,1,2\text{ ,} \\
\mathcal{N}_{u}(f_{\ell };\mathbb{S}^{2}) &=&\mathcal{N}_{u}\text{ , }Len(f_{\ell };\mathbb{S}^{2})=Len(0)\text{ ,}
\end{eqnarray*}%
so that $\mathcal{N}_{-\infty }$ is the total number of critical points and $%
\mathcal{L}_{2}(0)$ is the excursion area for $u=0.$ The correlation results
that we discussed so far can be summarized in the following table; here, we
denote by $u_{1},u_{2}\neq 0,1$ any two non-singular threshold values.

\begin{center}
The limiting value of $Corr^{2}(.,.),$ as $\ell \rightarrow \infty $

\begin{tabular}{llllllll}
& $\mathcal{L}_{j}(u_{1})$ & $\mathcal{L}_{j}(u_{2})$ & $Len(0)$ & $%
Len^{\ast }(u)$ & $\mathcal{L}_{2}(0)$ & $\mathcal{N}_{u}$ & $\mathcal{N}%
_{-\infty }$ \\
$\mathcal{L}_{j}(u_{1})$ & 1 & 1 & 0 & 0 & 0 & 1 & 0 \\
$\mathcal{L}_{j}(u_{2})$ & 1 & 1 & 0 & 0 & 0 & 1 & 0 \\
$Len(0)$ & 0 & 0 & 1 & 1 & 0 & 0 & 1 \\
$Len^{\ast }(u)$ & 0 & 0 & 1 & 1 & 0 & 0 & 1 \\
$\mathcal{L}_{2}(0)$ & 0 & 0 & 0 & 0 & 1 & 0 & 0 \\
$\mathcal{N}_{u}$ & 1 & 1 & 0 & 0 & 0 & 1 & 0 \\
$\mathcal{N}_{-\infty }$ & 0 & 0 & 1 & 1 & 0 & 0 & 1%
\end{tabular}%
$.$
\end{center}

\section{Eigenfunctions on Different Domains}

For brevity and simplicity's sake, this survey has focussed only on the
behaviour of random eigenfunctions on the sphere. Of course, as mentioned in
the Introduction this is just a special case of a much broader research
area, including for instance eigenfunctions on $\mathbb{R}^{d}$ and on the
standard flat torus $\mathbb{T}^{d}:=\mathbb{R}^d/\mathbb{Z}^d$. We do not even attempt to do justice to these
developments, but it is important to mention some of them which are
particularly close to the results we discussed for $\mathbb{S}^{2}$.

\subsection{Eigenfunctions on the Torus: Arithmetic Random Waves}

Eigenfunctions on the torus were first introduced in \cite{RW08} and have
then been studied by several other authors, see for instance \cite%
{Cam19,granville,KKW,Maff,MPRW2015,Rudnick,RudnickYesha} and the references
therein. In dimension $2$ these eigenfunctions (\emph{Arithmetic Random Waves%
}) are defined by the equations%
\begin{equation*}
\Delta _{\mathbb{T}^{2}}f_{n}+E_{n}f_{n}=0\text{ , }E_{n}=4\pi n\text{ , }%
n=a^{2}+b^{2}\text{ ,}
\end{equation*}%
for $a,b\in \mathbb{Z};$ the dimension of the $n$-th eigenspace is $\mathcal{%
N}_{n}:=Card\left\{ a,b\in \mathbb{Z}:a^{2}+b^{2}=n\right\} $, while the
expected value of nodal lengths is (\cite{RW08})%
\begin{equation*}
\mathbb{E}\left[ Len(f_{n};\mathbb{T}^{2})\right] =\frac{\sqrt{E_{n}}}{2%
\sqrt{2}}\text{ .}
\end{equation*}%
A major breakthrough was then obtained with the derivation of the variance
in \cite{KKW}. In this paper, the authors introduce a probability measure on
$\mathbb{S}^{1}$ defined by%
\begin{equation*}
\mu _{n}(.):=\frac{1}{\mathcal{%
N}_{n}}\sum_{a,b:a^{2}+b^{2}=n}\delta _{(a,b)}(\cdot)%
\text{ ,}
\end{equation*}%
$\delta _{(a,b)}(.)$ denoting the Dirac measure; its $k$-th order Fourier
coefficients are defined by $\widehat{\mu }_{n}(k):=\int_{\mathbb{S}%
^{1}}\exp (ik\theta )\mu _{n}(d\theta )$. In \cite{KKW} it is then shown
that the variance of nodal lengths has a non-universal behaviour and is
proportional to%
\begin{equation*}
Var\left\{ Len(f_{n};\mathbb{T}^{2})\right\} =\frac{1+\widehat{\mu }%
_{n}(4)^{4}}{512}\frac{E_{n}}{N_{n}^{2}}+o\left(\frac{E_{n}}{N_{n}^{2}}\right)\text{ ,
as }n\rightarrow \infty \text{ s.t. }\mathcal{%
N}_{n} \rightarrow \infty .
\end{equation*}%
It was later shown by \cite{MPRW2015} that the behaviour of $Len(\mathbb{T}%
^{2},f_{n})$ is dominated by its fourth-order chaos component, similarly to
what we observed above for random spherical harmonics (the result on the
torus was actually established earlier than the corresponding case for the
sphere). More precisely, we have that%
\begin{equation*}
Len(f_{n};\mathbb{T}^{2})-\mathbb{E}[Len(f_{n};\mathbb{T}^{2})]=\sum_{q=2}^{%
\infty }\text{Proj}\left[ \left. Len(f_{n};\mathbb{T}^{2})\right\vert 2q%
\right]
\end{equation*}%
\begin{equation*}
=\text{Proj}\left[ \left. Len(f_{n};\mathbb{T}^{2})\right\vert 4\right]
+o_{p}(\sqrt{Var\left\{ Len(f_{n};\mathbb{T}^{2})\right\} })\text{ ,}
\end{equation*}%
where Proj$\left[ \left. .\right\vert q\right] $ denotes projection on the $q
$-th order chaos. On the contrary of what we observed for the case of the
sphere, here it is not possible to express the fourth-order chaos as a
polynomial functional of the random eigenfunctions $\left\{ f_{n}\right\} $
alone. Moreover, the limiting distribution is non-Gaussian and
non-universal, i.e. it depends on the asymptotic behaviour of $%
\lim_{j\rightarrow \infty }\widehat{\mu }_{n_{j}}(4)$ which varies along
different subsequences $\left\{ n_{j}\right\} _{j=1,2,...}$(the attainable
measures for the weak-convergence of the sequences $\left\{ \mu
_{n_{j}}(.)\right\} _{n\in \mathbb{N}}$ have been investigated by \cite%
{KKW,KW}). Further results in this area include \cite{Cam19}, \cite{Not} for
arithmetic random waves in higher-dimension and \cite{WigYesha} for the
excursion area on subdomains of $\mathbb{T}^{2};$ as mentioned earlier, an
extension of Theorem \ref{2GKF} to the torus has been given in \cite{CMR}.
It should be noted that Arithmetic Random Waves can be viewed as an instance
of random trigonometric polynomials, whose zeroes have been studied, among
others, by \cite{Angst}, \cite{Bally}.

\subsection{The Euclidean Case: Berry's Random Waves}

Spherical harmonics on the sphere $\mathbb{S}^2$ are known to exhibit a scaling limit, i.e.
after a change of coordinates they converge locally to a Gaussian random
process on $\mathbb{R}^{2}$ which is isotropic, zero mean and has covariance
function%
\begin{equation*}
\mathbb{E}\left[ f(x)f(y)\right] =J_{0}(2\pi \left\Vert x-y\right\Vert ),%
\text{ }x,y\in \mathbb{R}^{2},\text{ }J_{0}(z):=\sum_{k=0}^{\infty }\frac{%
(-1)^{k}z^{2k}}{(k!)^{2}2^{2k}}\text{ ;}
\end{equation*}%
here $J_{0}(.)$ corresponds to the standard Bessel functions, for which the
following scaling asymptotics holds:%
\begin{equation*}
P_{\ell }(\cos \frac{\psi }{\ell })\rightarrow _{\ell \rightarrow \infty
}J_{0}(\psi )\text{ },\text{ }\psi \in \mathbb{R}\text{ }.
\end{equation*}
The behaviour of nodal lines $\mathcal{L}_{E}(f)=\left\{ x\in \mathbb{R}%
^{2}:f(x)=0,\left\Vert x\right\Vert <2\pi \sqrt{E}\right\} $ can then be
studied in the asymptotic regime $E\rightarrow \infty ;$ this is indeed the
physical setting under which Berry first investigated cancellation phenomena
in his pioneering paper \cite{Berry 2002}. The topology of nodal sets for
Berry's random waves was studied by \cite{nazarov}, \cite{Sarnak}, \cite{CH}
and others. Concerning nodal lengths, a (Quantitative) Central Limit Theorem
was established in \cite{npr}, where intersections of independent random
waves were also investigated; more recently, \cite{Vidotto21} proved a
result analogous to Theorem \ref{nodal}, namely that, as $E \rightarrow \infty$
\[
\mathcal{L}_{E}(f)-\mathbb{E}\left[\mathcal{L}_{E}(f)\right]
\]
\begin{equation}
=-\frac{1}{4}\frac{2\pi }{4!}\sqrt{\frac{E}{2}}%
\int_{\left\Vert x\right\Vert <2\pi \sqrt{E}}H_{4}(f(x))dx+o_{p}(\sqrt{%
Var\left\{ \mathcal{L}_{E}(f)\right\} })\text{ .}  \label{vidotto}
\end{equation}%
We expect that results analogous to (\ref{MRW}) and (\ref{vidotto}) will hold
for more general Riemannian waves on two-dimensional manifolds \cite{Zel};
extensions to random waves in $\mathbb{R}^{3}$ have been studied, among
others, by \cite{DalmaoEstradeLeon}, but in these higher-dimensional
settings it is no longer the case that nodal volumes are dominated by a
single chaotic component.

\subsection{Shrinking Domains}

As a final issue, we recall how some of the previous results can be extended
to shrinking subdomains of the torus and of the sphere. In this respect, a
surprising result was derived in \cite{BMW} concerning the asymptotic
behaviour of the nodal length on a suitably shrinking subdomain $%
B_{n}\subset \mathbb{T}^{2};$ indeed it was shown that, for density one
subsequences in $n$%
\begin{equation*}
\lim_{n\rightarrow \infty }Corr(Len(\mathbb{T}^{2},f_{n}),Len(\mathbb{T}%
^{2}\cap B_{n},f_{n}))=1\text{ ,}
\end{equation*}%
entailing that the behaviour of the nodal length on the whole torus is fully
determined by its behaviour on any shrinking disk $B_{n},$ provided the
radius of this disk is not smaller than $n^{-1/2+\varepsilon },$ some $%
\varepsilon >0.$ Of course, the asymptotic variance and distributions of the
Nodal Length in this shrinking domain is then immediately shown to be the
same as those for the full torus, up to a normalizing factor. Interestingly,
the same phenomenon does not occur on the sphere, where on the contrary it
was shown in \cite{Todino20} that%
\begin{equation*}
\lim_{\ell \rightarrow \infty }Corr(Len(\mathbb{S}^{2},f_{\ell }),Len(%
\mathbb{S}^{2}\cap B_{\ell },f_{\ell }))=0\text{ ,}
\end{equation*}%
so that the nodal length when evaluated on a shrinking subset $B_{\ell }$ of
the two-dimensional sphere is actually asymptotically independent from its
global value; in the same paper, it is indeed shown that (\ref{MRW})
generalizes to%
\begin{eqnarray}
&&Len(\mathbb{S}^{2}\cap B_{\ell },f_{\ell })-\mathbb{E}\left[ Len(\mathbb{S}%
^{2}\cap B_{\ell },f_{\ell })\right]  \label{todino} \\
&=&-\frac{1}{4}\sqrt{\frac{\lambda _{\ell }}{2}}\frac{1}{4!}h_{\ell
;4}(B_{\ell })+o_{p}(\sqrt{Var\left\{ h_{\ell ;4}(B_{\ell })\right\} })\text{
,}  \label{Todino2020}
\end{eqnarray}%
\begin{equation*}
h_{\ell ;4}(B_{\ell })=\int_{B_{\ell }}H_{4}(f_{\ell }(x))dx\text{ ;}
\end{equation*}%
from this characterization, a Central Limit Theorem follows easily along the
same lines that we discussed in the previous Section, see \cite{Todino20}
for more details and discussion.

\begin{acknowledgments}
I am grateful to Valentina Cammarota, Maurizia Rossi, Anna Paola Todino and an anonymous referee for a number of comments and suggestions on an earlier draft. This research was partly supported by the MIUR Departments of Excellence Program \emph{Math@Tov}.
\end{acknowledgments}

\end{document}